\documentclass[runningheads,a4paper]{llncs} 
\usepackage{url}
\usepackage{graphicx}
\usepackage{listings}
\usepackage{amsmath,amssymb,amsbsy,latexsym}
\newcommand{\nE}{\mathcal{E}}
\newcommand{\nF}{\mathcal{F}}
\newcommand{\nS}{\mathcal{S}}
\newcommand{\nD}{\mathcal{D}}
\newcommand{\nL}{\mathcal{L}}

\newcommand{\nP}{\widehat{\mathcal{L}}}
\newcommand{\dd}{{\rm{d}}}

\newcommand{\half}{\tfrac{1}{2}}

\newcommand\Order{\mathcal{O}}

\begin{document}

\lstset{language=matlab} 

%
%
\mainmatter             
\pagestyle{myheadings}  
\title{Symbolic Manipulation of Flows of Nonlinear Evolution Equations,
with Application in the Analysis of Split-Step Time Integrators}
\titlerunning{Symbolic Manipulation of Nonlinear Flows}  
%
\author{Winfried Auzinger \inst{1} \and Harald Hofst{\"a}tter \inst{1}
\and Othmar Koch \inst{2}}
\authorrunning{Winfried Auzinger et al.} 
\institute{Technische Universit{\"a}t Wien, Austria \\
\email{w.auzinger@tuwien.ac.at, hofi@harald-hofstaetter.at}, \\
\url{www.asc.tuwien.ac.at/~winfried},
\url{www.harald-hofstaetter.at}
\and
Universit{\"a}t Wien, Austria \\
\email{othmar@othmar-koch.org}, \\
\url{www.othmar-koch.org}
}

\maketitle              

\begin{abstract}
We describe a package realized in the Julia programming language which
performs symbolic manipulations applied to
nonlinear evolution equations, their flows, and commutators of such objects.
This tool was employed to perform contrived computations arising
in the analysis of the local error of operator splitting methods.
It enabled the proof of the convergence of the basic method and
of the asymptotical correctness of a defect-based error estimator.
The performance of our package is illustrated on several examples.
\end{abstract}
\section{Problem setting}
We are interested in the solution to nonlinear evolution equations
\begin{equation}\label{eq1}
\partial_t u(t) = A(u(t))+B(u(t)) = H(u(t)), \quad u(0) = u_0,
\end{equation}
on a Banach space $X$, where $A$ and $B$ are general nonlinear,
unbounded operators defined
on a 
subset $D\subset X$, the solution is denoted by $\mathcal{E}_H(t,u_0)$,
and analogously for the two sub-flows associated with $A$ and $B$. The structure of the vector fields often suggests
to employ additive splitting methods to separately propagate the two subproblems
defined by $A$ and $B$,
\begin{equation}\label{splitting1}
u(t_1)\approx u_{1} := \mathcal{S}(h,u_0) = \nE_B(b_s h,\cdot) \circ \nE_A(a_s h,\cdot) \circ \ldots \circ
                              \nE_B(b_1 h,\cdot) \circ \nE_A(a_1 h,u_0),
\end{equation}
where the coefficients $ a_j,b_j,\, j=1 \ldots s $ are determined according to the requirement
that a prescribed order of consistency is obtained~\cite{haireretal02}.

Both in the a priori error analysis and for a posteriori error estimation,
a defect-based approach has been introduced
in~\cite{auzingeretal13b},
which serves both to derive theoretical error bounds and as a basis for adaptive step-size
selection.

In the analysis of this error estimate, extensive symbolic manipulation of the
flows defined by the operators in~\eqref{eq1}, their Fr\'{e}chet derivatives and
arising commutators is indispensable. As such calculations imply a high
effort of formula manipulation and are highly error prone, we have implemented
an automatic tool in the Julia programming language to verify the calculation.
Julia is a high-level, high-performance dynamic programming language
for technical computing, see~\cite{julia}.

By defining a suitable set
of substitution rules, the algorithm can check the equivalence of expressions
built from the objects mentioned. On this basis, it was possible to
ascertain the correctness of all steps in the proofs given in
\cite{auzingeretal13b}.

\section{Local error, defect, and error estimator}\label{section2}

We describe the background arising from the application of splitting methods~\eqref{splitting1}
to the solution of evolution equations~\eqref{eq1}, see~\cite{auzingeretal13b}. The defect of the splitting
approximation is defined as
\begin{equation*} 
\nD(t,u) = \nS^{(1)}(t,u) = \partial_t \nS(t,u) - H(\nS(t,u)),
\end{equation*}
while the local error is given by
\begin{equation} \label{locerr1}
 \nL(t,u) = \nS(t,u)-\nE_H(t,u) = \int_0^t\,\nF(\tau,t,u)\,\dd\tau,
\end{equation}
with
\begin{equation*} 
\nF(\tau,t,u)=\partial_2\nE_H(t-\tau,\nS(\tau,u))\cdot\,\nS^{(1)}(\tau,u).
\end{equation*}
\subsection*{The Lie-Trotter method}
We illustrate our analysis of the local error for the simplest
Lie-Trotter splitting method, 
\begin{equation*}
\nS(t,u_0)=\nE_B(t,\nE_A(t,u_0)),
\end{equation*}
see~\cite{auzingeretal13b}.
This involves some nontrivial crucial identities, 
namely~\eqref{check1}--\eqref{check5} below, which will be
verified using our Julia package in Section~\ref{sec:verify}.

Our aim is to show
\begin{equation*}
\nS^{(1)}(t,u)=\Order(t), 
\end{equation*}
and thus
\begin{equation*}
\nL(t,u) = \Order(t^2).
\end{equation*}
$ \nS^{(1)}(t,u) $ can be represented in the form
\begin{subequations}\label{check1}
\begin{align}
 \nS^{(1)}(t,u) &= \partial_t\,\nS(t,u) -H(\nS(t,u))
  = \tilde\nS^{(1)}(t,\nE_A(t,u)),
  \label{check_S1}
\end{align}
with
\begin{equation}
\label{S1schlange}
\tilde\nS^{(1)}(t,v) = \partial_2\nE_B(t,v)\cdot A(v) - A(\nE_B(t,v)).
\end{equation}
\end{subequations}
$\tilde\nS^{(1)}(t,v)$ satisfies
\begin{equation}
\label{check2}
\begin{aligned}
\partial_t\tilde\nS^{(1)}(t,v) &=
B'(\nE_B(t,v))\cdot\tilde\nS^{(1)}(t,v)+[B,A](\nE_B(t,v)), \\
\tilde\nS^{(1)}(0,v) &= 0,
\end{aligned}
\end{equation}
where
\begin{equation*}
[B,A](u) = B'(u)A(u) - A'(u)B(u)
\end{equation*}
denotes the commutator of the two vector fields.

%
From
\begin{equation*}
\partial_t \nE_F(t,u) = F(\nE_F(t,u))
~~\Rightarrow~~
\partial_t \partial_{2}\nE_{F}(t,u)\cdot v
  =F'(\nE_{F}(t,u))\cdot \partial_{2}\nE_{F}(t,u)\cdot v
\end{equation*}
it follows 
that $\partial_2\nE_F(t,u)$ is a fundamental system of the linear differential equation
\begin{equation*}
  \partial_tX(t,u)=
F'(\nE_F(t,u))\cdot X(t,u)
\end{equation*}
which  satisfies
\begin{equation*}
  \partial_2\nE_F(t,u)^{-1}
  = \partial_2\nE_F(-t,\nE_F(t,u)).
\end{equation*}
The solution of an inhomogenous system like (\ref{check2}) of the form
\begin{align*}
&\partial_tX(t,u)=
F'(\nE_F(t,u))\cdot X(t,u) + R(t,u),\\
&X(0,u)=X_0(u)
\end{align*}
can be represented by the variation of constant formula,
\begin{equation*}
  X(t,u) = \partial_2\nE_F(t,u)\cdot\left(X_0(u)+\int_0^t
  \partial_2\nE_F(-\tau,\nE_F(\tau,u))\cdot R(\tau,u)\,\dd\tau\right).
\end{equation*}
Hence, the term $ \tilde\nS^{(1)}(t,v) $ defined in~\eqref{S1schlange} satisfies
\begin{equation*}
\tilde\nS^{(1)}(t,v) =  \partial_2\nE_B(t,v)\cdot
  \int_0^t\partial_2\nE_B(-\tau,\nE_B(\tau,v) )\cdot
  [B,A](\nE_B(\tau,v))\,\dd\tau.
\end{equation*}
From this integral representation it follows
\begin{equation*}
\nD(t,u) = \tilde\nS^{(1)}(t,\nE_A(t,u)) = \Order(t),
\end{equation*}
and
\begin{equation*}
\nL(t,u) = \int_0^t\partial_2\nE_H(t-\tau,\nS(\tau,u))\cdot
\nD(\tau,u)\,\dd\tau
= \Order(t^2).
\end{equation*}

As a basis for adaptive time-stepping, we define an
a~posteriori local error estimator
by numerical evaluation of the integral representation~\eqref{locerr1} of the
local error by the trapezoidal rule, yielding
\begin{equation*}
 \nP(t,u) = \half\,t\,\nF(t,t,u) = \half\,t\,\nD(t,u)
 =\half\,t\,\nS^{(1)}(t,u).
\end{equation*}
To analyze the deviation of this error estimator from the exact error,
we use the Peano representation
\begin{equation*}
  \nP(t,u)- \nL(t,u) =
  \int_0^t K_1(\tau,t)\,\partial_\tau\nF(\tau,t,u)\,\dd\tau,
\end{equation*}
with the kernel
\begin{equation*}
  K_1(\tau,t) = \tau-\half\,t = \Order(t).
\end{equation*}
To infer asymptotical correctness of the error estimator, we wish to show that
\begin{equation*}
 \nP(t,u)- \nL(t,u) = \Order(t^3).
\end{equation*}
To this end, we compute
\begin{align}
\partial_\tau\nF(\tau,t,u)
   &=\partial_2\nE_H(t-\tau,\nS(\tau,u))\cdot\,\nS^{(2)}(\tau,u)
   \nonumber\\
   & \label{check2a} \quad {} + \partial_2^2\nE_H(t-\tau,\nS(\tau,u))(\nS^{(1)}(\tau,u),\nS^{(1)}(\tau,u)) \\
   &=\partial_2\nE_H(t-\tau,\nS(\tau,u))\cdot\,\nS^{(2)}(\tau,u) + \Order(t),\nonumber
\end{align}
where
\begin{subequations}
\label{check3}
\begin{align}
\nS^{(2)}(t,u)&=\partial_t\,\nS^{(1)}(t,u)
-H'(\nS(t,u))\cdot \nS^{(1)}(t,u)\nonumber\\
&=\tilde\nS^{(2)}(t,\nE_A(t,u)),
\end{align}
with
\begin{equation}
\tilde\nS^{(2)}(t,v) =
\partial_2\tilde\nS^{(1)}(t,v)\cdot A(v)
-A'(\nE_B(t,v))\cdot\tilde\nS^{(1)}(t,v)+[B,A](\nE_B(t,v)).
\end{equation}
\end{subequations}
$\tilde\nS^{(2)}(t,v)$ satisfies
\begin{align}
\partial_t\tilde\nS^{(2)}(t,v)
&= B'(\nE_B(t,v))\cdot\tilde\nS^{(2)}(t,v)\nonumber\\
& \quad {} +B''(\nE_B(t,v))(\tilde\nS^{(1)}(t,v),\tilde\nS^{(1)}(t,v))\nonumber\\
& \quad {} -[B,[B,A]](\nE_B(t,v)) -[A,[B,A]](\nE_B(t,v)) \label{check4} \\
& \quad {} +2[B,A]'(\nE_B(t,v))\cdot\tilde\nS^{(1)}(t,v), \nonumber \\
\tilde \nS^{(2)}(0,v) &= [B,A](v).\nonumber
\end{align}
This implies the integral representation
\begin{align*}
\tilde\nS^{(2)}(t,v) &=
\partial_2\nE_B(t,v)\cdot [B,A](v) +
 \partial_2\nE_B(t,v)\cdot
  \int_0^t\partial_2\nE_B(-\tau,\nE_B(\tau,v) )\,\cdot\\
 & \quad \Big(
 B''(\nE_B(\tau,v))(\tilde\nS^{(1)}(\tau,v),\tilde\nS^{(1)}(\tau,v))\\
& \quad~ {} - [B,[B,A]](\nE_B(\tau,v)) -[A,[B,A]](\nE_B(\tau,v))\\
& \quad~ {} + 2[B,A]'(\nE_B(\tau,v))\cdot\tilde\nS^{(1)}(\tau,v) \Big)\,\dd\tau.
\end{align*}
Thus,
\begin{equation*}
\nS^{(2)}(\tau,u) = \partial_2\nE_B(\tau,\nE_A(\tau,u))\cdot [B,A](\nE_A(\tau,u)) + \Order(t),
\end{equation*}
and altogether
{
\begin{align*}
&\nP(t,u)- \nL(t,u) =
  \int_0^t K_1(\tau,t)\,\partial_\tau\nF(\tau,t,u)\,\dd\tau \\
&=
  \int_0^t K_1(\tau,t)\,\,\cdot \\
& \quad \cdot\,\partial_2\nE_H(t-\tau,\nS(\tau,u))\cdot
  \partial_2\nE_B(\tau,\nE_A(\tau,u))\,\cdot\,
                              [B,A](\nE_A(\tau,u))\,\dd\tau + \Order(t^3) \\
&=  \int_0^t K_2(\tau,t)\,\,\cdot \\
& \quad \cdot\,\partial_\tau\Big(
  \partial_2\,\nE_H(t-\tau,\nS(\tau,u))\,\cdot\,
  \partial_2\,\nE_B(\tau,\nE_A(\tau,u))\cdot [B,A](\nE_A(\tau,u))
  \Big)\,\dd\tau + \Order(t^3),
\end{align*}
}
where $K_2(\tau,t) =\half\tau(t-\tau) = \Order(t^2)$ by partial integration.
Here,
\begin{subequations}\label{check5}
\begin{equation}
  \partial_\tau\Big(
  \partial_2\nE_H(t-\tau,\nS(\tau,u))\cdot
  \partial_2\nE_B(\tau,\nE_A(\tau,u))\cdot [B,A](\nE_A(\tau,u))
  \Big) = \Order(1),
\end{equation}
because this derivative can be expressed as
\begin{align}
&\Big[\partial_2\nE_H(t-\tau,\nE_B(\tau,v))\cdot\partial_2\nE_B(\tau,v)\cdot[[B,A],A](v)\nonumber\\
& \qquad {} + \partial_2\nE_H(t-\tau,\nE_B(\tau,v))\cdot\partial_2\tilde\nS^{(1)}(\tau,v)\cdot[B,A](v)\nonumber\\
& \qquad {} + \partial_2^2\nE_H(t-\tau,\nE_B(\tau,v))\big(\tilde\nS^{(1)}(\tau,v), \partial_2\nE_B(\tau,v)\cdot[B,A](v)\big)
\Big]_{v=\nE_A(t,u)}.
\end{align}
\end{subequations}

\subsection*{The Strang splitting method}
For the Strang splitting method
\begin{equation*}
\nS(t,u_0)=\nE_B(\tfrac{t}{2},\nE_A(t,\nE_B(\tfrac{t}{2},u_0))),
\end{equation*}
our aim is to show
\begin{equation*}
\nS^{(1)}(t,u)=\Order(t^2),
\end{equation*}
and thus
\begin{equation*}
 \nL(t,u) = \Order(t^3).
\end{equation*}
In this case the necessary manipulations become significantly more complex
than for the Lie-Trotter case, in particular concerning the analysis of
a defect-based a~posteriori error estimator.
The analysis is based on multiple application
of variation of constant formulas in a general nonlinear setting.
This was the main motivation for the development of our computational tool;
in particular, the theoretical results from~\cite{auzingeretal13b},
which are not repeated here, have been verified using this tool.

\section{The Julia package {\tt Flows.jl}}

The package described in the following is available from~\cite{Flowsjl}.
It consists of approximately 1000 lines of Julia code and is essentially self-contained,
relying only on the Julia standard library but not on additional packages.
A predecessor written in Perl has been used for the preparation
of~\cite{auzingeretal13b}.
In a Julia notebook, the package is initialized as follows:
\begin{itemize}
\item[] {\tt In [1]: using Flows}
\end{itemize}

\begin{figure}[!t]
\begin{center}
\vspace{-12mm}
\includegraphics[width=0.71\textwidth,angle=-90]{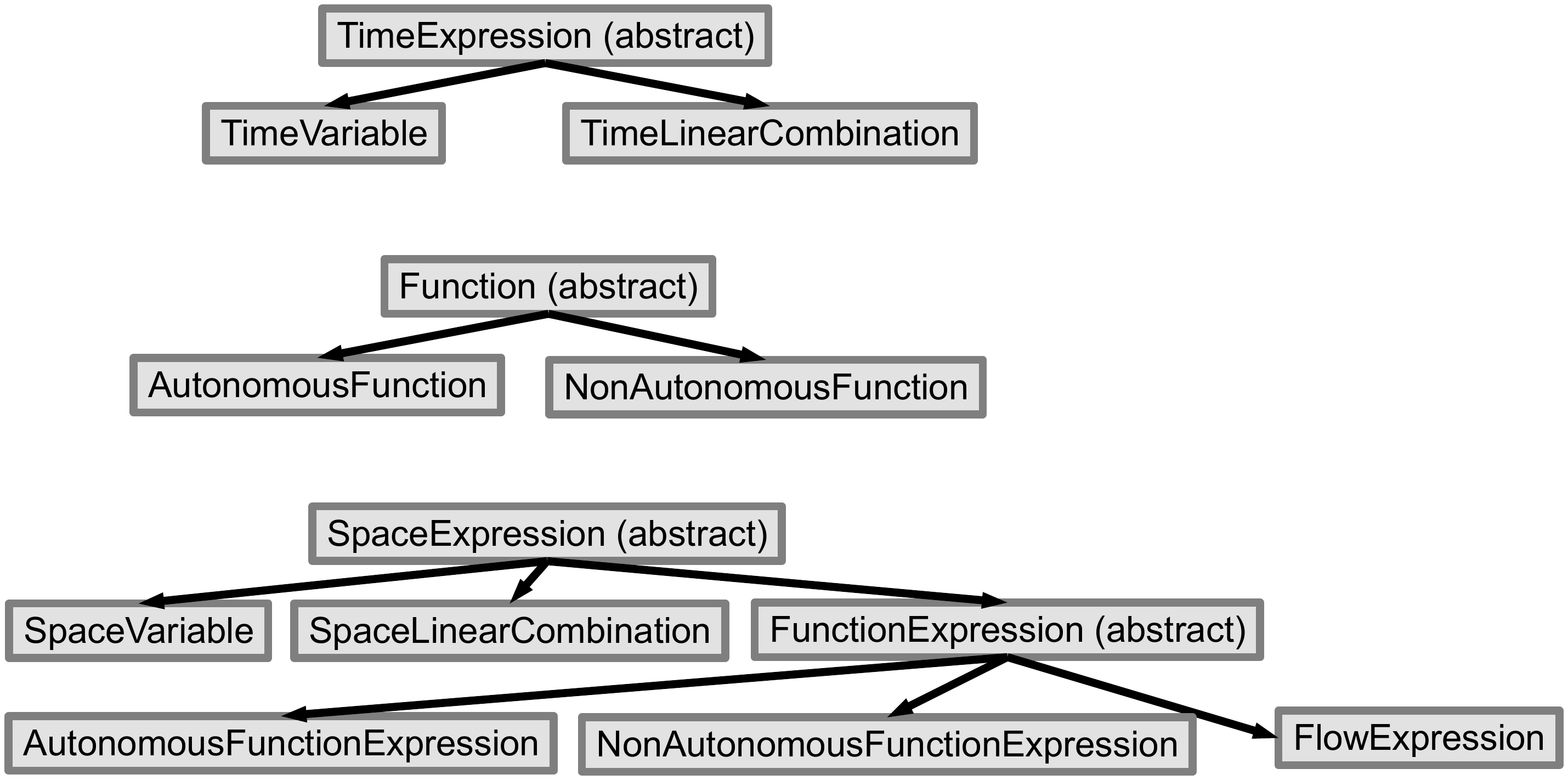}
\vspace{-12mm}
\caption{Data type hierarchy of {\tt Flows.jl}.\label{fig:hierarchy1}}
\end{center}
\end{figure}

\subsection*{Data types}
The data types in the package {\tt Flows.jl} and their hierarchical dependence
are illustrated in~Figure~\ref{fig:hierarchy1}.

\medskip\noindent
Objects of the abstract type {\tt TimeExpression}
represent a first argument $ t $ in a flow expression like $ \nE_H(t,u) $.

\medskip\noindent
Objects of type {\tt TimeVariable} are generated as follows:
\begin{itemize}
\item[] {\tt In [2]: @t\_vars t s r}
\item[] {\tt Out[2]: (t,s,r)}
\end{itemize}
Objects of type {\tt TimeLinearCombination}:
\begin{itemize}
\item[] {\tt In [3]: ex = t - 2s + 3r}
\item[] {\tt Out[3]: $3r + t - 2s$}
\end{itemize}

\smallskip\noindent
Similarly, objects of the abstract type {\tt SpaceExpression}
represent a second argument $ u $ in a flow expression like $ \nE_H(t,u) $.

\smallskip\noindent
Objects of type {\tt SpaceVariables}:
\begin{itemize}
\item[] {\tt In [4]: @x\_vars u v w}
\item[] {\tt Out[4]: (u,v,w)}
\end{itemize}

\smallskip\noindent
In order to construct objects of type
{\tt AutonomousFunctionExpression} or {\tt Flow\-Expression}
like $ A(u) $ or $ \nE_A(t,u) $ we first need to declare a symbol for $ A $
of type {\tt AutonomousFunction}.\footnote{Likewise for objects of type {\tt NonAutonomousFunction}.}
\begin{itemize}
\item[] {\tt In [5]: @funs A}
\item[] {\tt Out[5]: (A,)}
\end{itemize}
Now we can generate objects of types {\tt AutonomousFunctionExpression} and {\tt FlowExpression}:
\begin{itemize}
\item[] {\tt In [6]: ex1 = A(u)}
\item[] {\tt Out[6]: $A(u)$}
\\[-3mm]
\item[] {\tt In [7]: ex2 = E(A,t,u)}
\item[] {\tt Out[7]: $\nE_A(t,u)$}
\end{itemize}
Additional arguments in such expressions represent arguments for Fr{\'e}chet derivatives
with respect to a {\tt SpaceVariable} {\tt u}. Note that the order of the derivative
is implicitly determined from the number of arguments.
\begin{itemize}
\item[] {\tt In [8]: ex1 = A(u,v,w)}
\item[] {\tt Out[8]: $A''(u)(v,w)$}
\\[-3mm]
\item[] {\tt In [9]: ex2 = E(A,t,u,v)}
\item[] {\tt Out[9]: $\partial_2\nE_A(t,u) \cdot v$}
\end{itemize}
Objects of type {\tt SpaceLinearCombination} can be built from such expressions:
\begin{itemize}
\item[] {\tt In [10]: ex = -2E(A,t,u,v,w)+2u+A(v,w)}
\item[] {\tt Out[10]: $-2\partial_2^2\nE_A(t,u)(v,w)+2u+A'(v) \cdot w $}
\end{itemize}

\subsection*{Methods}

\subsubsection*{{\tt differential}}
$~$\\[1.5mm]
This method generates the Fr{\'e}chet derivative of an expression with respect to a {\tt SpaceVariable}
applied to an expression.
\begin{itemize}
\item[] {\tt In [11]: ex = A(B(u)) + E(A,t,u,v)}
\item[] {\tt Out[11]: $A(B(u)) + \partial_2 \nE_A(t,u) \cdot v$}
\\[-3mm]
\item[] {\tt In [12]: differential(ex,u,B(w))}
\item[] {\tt Out[12]: $A'(B(u)) \cdot B'(u) \cdot B(w) + \partial_2^2 \nE_A(t,u)(v,B(w))$}
\end{itemize}

\subsubsection*{{\tt t\_derivative}}
$~$\\[1.5mm]
This method generates the derivative of an expression with respect to a {\tt Time\-Variable}.
\begin{itemize}
\item[] {\tt In [13]: ex = E(A,t-2s,u+E(B,s,v))}
\item[] {\tt Out[13]: $\nE_A(t-2s,u+\nE_B(s,v))$}
\\[-3mm]
\item[] {\tt In [14]: t\_derivative(ex,s)}
\item[] {\tt Out[14]: $\partial_2 \nE_A(t-2s,u+\nE_B(s,v)) \cdot B(\nE_B(s,v)) - 2A(\nE_A(t-2s,u+\nE_B(s,v))$}
\end{itemize}

\subsubsection*{{\tt expand}}
$~$\\[1.5mm]
A (higher-order) Fr{\'e}chet derivative is a (multi-)linear map.
The method {\tt expand} transforms the application of such a (multi-)linear map
to a linear combination of expressions into the corresponding linear
combination of (multi-)linear maps.
\begin{itemize}
\item[] {\tt In [15]: ex1 = E(A,t,u,2v+3w)}
\item[] {\tt Out[15]: $\partial_2 \nE_A(t,u) \cdot (2v+3w)$}
\\[-3mm]
\item[] {\tt In [16]: expand(ex1)}
\item[] {\tt Out[16]: $2\partial_2 \nE_A(t,u) \cdot v + 3\partial_2 \nE_A(t,u) \cdot w $}
\\[-3mm]
\item[] {\tt In [17]: ex2 = A''(u,v+w,v+w)}
\item[] {\tt Out[17]: $2A''(u)(v+w,v+w)$}
\\[-3mm]
\item[] {\tt In [18]: expand(ex2)}
\item[] {\tt Out[18]: $2A''(u)(v,w) + A''(u)(v,v) + A''(u)(w,w)$}
\end{itemize}

\subsubsection*{{\tt substitute}}
$~$\\[1.5mm]
Different variants of substitutions are implemented.
The most sophisticated one is substitution of an object of type
{\tt Function} by an expression.
For example, this allows to define the double commutator $[A,[A,B]](u)$ by substituting
$B$ in 
$[A,B](u)$ by $[A,B](u)$:
\begin{itemize}
\item[] {\tt In [19]: C\_AB = A(u,B(u))-B(u,A(u))}
\item[] {\tt Out[19]: $A'(u) \cdot B(u) - B'(u) \cdot A(u)$}
\\[-3mm]
\item[] {\tt In [20]: substitute(C\_AB,B,C\_AB,u)}
\item[] {\tt Out[20]: $-A''(u)(B(u),A(u))+B'(u)\cdot A'(u)\cdot A(u)+B''(u)(A(u),A(u))$}

\qquad\quad\, \,  ~ $-A'(u)\cdot B'(u)\cdot A(u)+A'(u)\cdot (-B'(u)\cdot A(u)+A'(u)\cdot B(u)) $
\end{itemize}

\subsubsection*{{\tt commutator}}
$~$\\[1.5mm]
This method generates expressions involving commutators $[A,B]$ and
double commutators $[A,[B,C]]$.
\begin{itemize}
\item[] {\tt In [21]: ex1 = commutator(A,B,u) }
\item[] {\tt Out[21]: $A'(u) \cdot B(u) - B'(u) \cdot A(u)$}
\pagebreak
\item[] {\tt In [22]: ex2 = commutator(A,B,C,u)}\nopagebreak\nopagebreak
\item[]\nopagebreak {\tt Out[22]: $C'(u)\cdot B'(u)\cdot A(u)+C''(u)(B(u),A(u))-B'(u)\cdot C'(u)\cdot A(u)$}

\qquad\qquad\ ${}-B''(u)(C(u),A(u))+A'(u)\cdot (B'(u)\cdot C(u)-C'(u)\cdot B(u))$
\end{itemize}
Example: We verify the Jacobi identity $[A,[B,C]]+[B,[C,A]]+[C,[A,B]]=0$:
\begin{itemize}
\item[] {\tt In [23]: ex3 = expand(commutator(A,B,C,u)+commutator(B,C,A,u)

        \qquad\qquad\qquad \, \,  {\tt +commutator(C,A,B,u}))}
\item[] {\tt Out[23]:} $ 0 $
\end{itemize}

\subsubsection*{{\tt FE2DEF}, {\tt DEF2FE}}
$~$\\[1.5mm]
These methods substitute expressions according to the fundamental identity
\begin{equation*}
A(\nE_{A}(t,u))=\partial_{2}\nE_{A}(t,u)\cdot A(u),
\end{equation*}
see \cite{auzingeretal13b}.
\begin{itemize}
\item[] {\tt In [24]: ex1 = A(E(A,t,u))}
\item[] {\tt Out[24]: $A(\mathcal{E}_{A}(t,u))$}
\\[-3mm]
\item[] {\tt In [25]: ex2 = FE2DEF(ex1)}
\item[] {\tt Out[25]: $\partial_{2}\mathcal{E}_{A}(t,u)\cdot A(u)$}
\\[-3mm]
\item[] {\tt In [26]: DEF2FE(ex2)}
\item[] {\tt Out[26]: $A(\mathcal{E}_{A}(t,u))$}
\end{itemize}

\subsubsection*{{\tt reduce\_order}}
$~$\\[1.5mm]
This method constitutes the essential manipulation needed for the verification of the
identities~\eqref{check1}--\eqref{check5}
in Section~\ref{section2}.
By repeated differentiation of the fundamental  identity
\begin{equation*}
A(\nE_{A}(t,u))-\partial_{2}\nE_{A}(t,u)\cdot A(u)=0
\end{equation*}
we obtain
\begin{eqnarray*}
&&A'(\nE_{A}(t,u))\cdot \partial_{2}\nE_{A}(t,u)\cdot v-\partial_{2}^{2}\nE_{A}(t,u)(A(u),v)-\partial_{2}\nE_{A}(t,u)\cdot A'(u)\cdot v=0,\\[1mm]
&&A''(\nE_{A}(t,u))(\partial_{2}\nE_{A}(t,u)\cdot v,\partial_{2}\nE_{A}(t,u)\cdot w)
+A'(\nE_{A}(t,u))\cdot \partial_{2}^{2}\nE_{A}(t,u)(v,w)\\
&&\qquad-\partial_{2}^{3}\nE_{A}(t,u)(A(u),v,w)
-\partial_{2}^{2}\nE_{A}(t,u)(A'(u)\cdot v,w)\\
&&\qquad-\partial_{2}^{2}\nE_{A}(t,u)(A'(u)\cdot w,v)
-\partial_{2}\nE_{A}(t,u)\cdot A''(u)(v,w)=0,
\end{eqnarray*}
and so on.
The method {\tt reduce\_order} transforms expressions of the form of the
highest order derivative in these identities by means of these very identities:
\begin{itemize}
\item[] {\tt In [27]: ex1 = E(A,t,u,A(u))}
\item[] {\tt Out[27]: $\partial_{2}\nE_{A}(t,u)\cdot A(u)$}
\\[-3mm]
\item[] {\tt In [28]: reduce\_order(ex1)}
\item[] {\tt Out[28]: $A(\nE_{A}(t,u))$}
\pagebreak
\item[] {\tt In [29]: ex2 = E(A,t,u,A(u),v)}
\item[] {\tt Out[29]: $\partial_{2}^{2}\nE_{A}(t,u)(A(u),v)$}
\\[-3mm]
\item[] {\tt In [30]: reduce\_order(ex2)}
\item[] {\tt Out[30]: $A'(\nE_{A}(t,u))\cdot \partial_{2}\nE_{A}(t,u)\cdot v
-\partial_{2}\nE_{A}(t,u)\cdot A'(u)\cdot v$}
\\[-3mm]
\item[] {\tt In [31]: ex3 = E(A,t,u,A(u),v,w)}
\item[] {\tt Out[31]: $\partial_{2}^{3}\nE_{A}(t,u)(A(u),v,w)$}
\\[-3mm]
\item[] {\tt In [32]: reduce\_order(ex3)}
\item[] {\tt Out[32]:}
\vspace{-6.25mm}
    \begin{align*}
%
 &\quad\, -\partial_{2}^{2}\mathcal{E}_{A}(t,u)(A'(u)\cdot w,v)-\partial_{2}^{2}\mathcal{E}_{A}(t,u)(A'(u)\cdot v,w)\\[-0.5\jot]
 &\quad\,-\partial_{2}\mathcal{E}_{A}(t,u)
 \cdot A''(u)(v,w)
 +A'(\mathcal{E}_{A}(t,u))\cdot \partial_{2}^{2}\mathcal{E}_{A}(t,u)(v,w)\\[-0.5\jot]
 &\quad\, +A''(\mathcal{E}_{A}(t,u))(\partial_{2}\mathcal{E}_{A}(t,u)\cdot v,\partial_{2}\mathcal{E}_{A}(t,u)\cdot w)   
    \end{align*}
\\[-11mm]
\end{itemize}
Similarly for higher derivatives of analogous form.

\section{Verification of crucial identities} \label{sec:verify}
We describe a Julia notebook which  implements the verification of the 
identities (\ref{check1})--(\ref{check5}) of Section~\ref{section2}.

\subsection*{Check~\eqref{check1}}

We verify identity~\eqref{check1},
\begin{equation*}
 \partial_t\nS(t,u) - H(\nS(t,u)) =\Big[\partial_2\nE_B(t,v)\cdot A(v) - A(\nE_B(t,v))
 \Big]_{v=\nE_A(t,u)}.
\end{equation*}

\begin{itemize}
\item[] {\tt In [{1}]: using Flows}
\\[-3mm]
\item[] {\tt In [{2}]: {@}{t\_vars} t  }
\item[] {\tt Out[{2}]: (t,)  }
\\[-3mm]
\item[] {\tt In [{3}]: {@}{x\_vars} u v}
\item[] {\tt Out[{3}]: (u,v)  }
\\[-3mm]
\item[] {\tt In [{4}]: {@}{funs} A B }
\item[] {\tt Out[{4}]: (A,B) }
\\[-3mm]
\item[] {\tt In [{5}]: {E\_Atu} {=} {E}{(}{A}{,}{t}{,}{u}{)}}
\item[] {\tt Out[{5}]: $\mathcal{E}_{A}(t,u)$}
\\[-3mm]
\item[] {\tt In [{6}]: {E\_Btv} {=} {E}{(}{B}{,}{t}{,}{v}{)}}
\item[] {\tt Out[{6}]: $\mathcal{E}_{B}(t,v)$}
\\[-3mm]
\item[] {\tt In [{7}]: {Stu} {=} {E}{(}{B}{,}{t}{,}{E}{(}{A}{,}{t}{,}{u}{)}{)}}
\item[] {\tt Out[{7}]: $\mathcal{E}_{B}(t,\mathcal{E}_{A}(t,u))$}
\\[-3mm]
\item[] {\tt In [{8}]: {S1tv} {=} {differential}{(}{E}{(}{B}{,}{t}{,}{v}{)}{,}{v}{,}{A}{(}{v}{)}{)}{-}{A}{(}{E}{(}{B}{,}{t}{,}{v}{)}{)}}
\item[] {\tt Out[{8}]: $\partial_{2}\mathcal{E}_{B}(t,v)\cdot A(v)-A(\mathcal{E}_{B}(t,v))$}
\\[-3mm]
\item[] {\tt In [{9}]: {S1tu} {=} {substitute}{(}{S1tv}{,}{v}{,}{E}{(}{A}{,}{t}{,}{u}{)}{)}}

        \qquad\quad\, \ \ {\tt {ex1} {=} {S1tu}}
\item[] {\tt Out[{9}]:}
    \(-A(\mathcal{E}_{B}(t,\mathcal{E}_{A}(t,u)))+\partial_{2}\mathcal{E}_{B}(t,\mathcal{E}_{A}(t,u))\cdot A(\mathcal{E}_{A}(t,u))\)
\\[-3mm]
\item[] {\tt In [{10}]: {ex2} {=} {t\_derivative}{(}{Stu}{,}{t}{)}{-}{(}{A}{(}{Stu}{)}{+}{B}{(}{Stu}{)}{)}}
\item[] {\tt Out[{10}]: $-A(\mathcal{E}_{B}(t,\mathcal{E}_{A}(t,u)))
+\partial_{2}\mathcal{E}_{B}(t,\mathcal{E}_{A}(t,u))\cdot A(\mathcal{E}_{A}(t,u))$}
\\[-3mm]
\item[] {\tt In [{11}]: {ex1}{-}{ex2}}
\item[] {\tt Out[{11}]: $0$}
\end{itemize}

\subsection*{Check~\eqref{check2}}
We verify identity~\eqref{check2},
\begin{equation*}
\partial_t\tilde\nS^{(1)}(t,v) =
B'(\nE_B(t,v))\cdot\tilde\nS^{(1)}(t,v)+[B,A](\nE_B(t,v)).
\end{equation*}

\begin{itemize}

\item[] {\tt In [{12}]: {ex1} {=} {B}{(}{E}{(}{B}{,}{t}{,}{v}{)}{,}{S1tv}{)}{+}{commutator}{(}{B}{,}{A}{,}{E}{(}{B}{,}{t}{,}{v}{)}{)}}
\item[] {\tt Out[12]: $-A'(\mathcal{E}_{B}(t,v))\cdot B(\mathcal{E}_{B}(t,v))+B'(\mathcal{E}_{B}(t,v))\cdot (\partial_{2}\mathcal{E}_{B}(t,v)\cdot A(v)$}
    
    \qquad\qquad \ \ \(-A(\mathcal{E}_{B}(t,v)))+B'(\mathcal{E}_{B}(t,v))\cdot A(\mathcal{E}_{B}(t,v))\)
\\[-3mm]
\item[] {\tt In [{13}]: {ex2} {=} {t\_derivative}{(}{S1tv}{,}{t}{)}}
\item[] {\tt Out[{13}]: $-A'(\mathcal{E}_{B}(t,v))\cdot B(\mathcal{E}_{B}(t,v))+B'(\mathcal{E}_{B}(t,v))\cdot \partial_{2}\mathcal{E}_{B}(t,v)\cdot A(v)$}
\\[-3mm]
\item[] {\tt In [{14}]: {reduce\_order}{(}{expand}{(}{ex1}{-}{ex2}{)}{)}}
\item[] {\tt Out[{14}]: $0$}
\end{itemize}

\subsection*{Check~\eqref{check2a}}
We verify identity~\eqref{check2a},
\begin{align*}
&\partial_\tau\big(\partial_2\nE_H(t-\tau,\nS(\tau,u))\cdot\,\nS^{(1)}(\tau,u)\big)\\
&=\partial_2\nE_H(t-\tau,\nS(\tau,u))\cdot\,\nS^{(2)}(\tau,u)
   \nonumber
    + \partial_2^2\nE_H(t-\tau,\nS(\tau,u))(\nS^{(1)}(\tau,u),\nS^{(1)}(\tau,u)).
\end{align*}

\begin{itemize}

\item[] {\tt In [{15}]: {@}{nonautonomous\_funs} {S}}
\item[] {\tt Out[{15}]: (S,)}
\\[-3mm]
\item[] {\tt In [{16}]: {@}{funs} {H}}
\item[] {\tt Out[{16}]: (H,)}
\\[-3mm]
\item[] {\tt In [{17}]: S1tu = t\_derivative(S(t,u),t)-H(S(t,u)) }
\item[] {\tt Out[{17}]: $-H(S(t,u))+\partial_{1}S(t,u)$}
\\[-3mm]
\item[] {\tt In [{18}]: S2tu = t\_derivative(S1tu,t)-H(S(t,u),S1tu)}
\item[] {\tt Out[{18}]: $-H'(S(t,u))\cdot \partial_{1}S(t,u)+\partial_{1}^{2}S(t,u)-H'(S(t,u))\cdot (-H(S(t,u))$}

    \qquad\qquad \ \ $+\partial_{1}S(t,u))$
\\[-3mm]
\item[] {\tt In [{19}]: {@}{t\_vars} {T}}
\item[] {\tt Out[{19}]: (T,)}
\end{itemize}
In the following, a trailing semicolon in the input suppresses the display of the corresponding output.
\begin{itemize}
\item[] {\tt In [{20}]: ex1 = E(H,T-t,S(t,u),S2tu)+E(H,T-t,S(t,u),S1tu,S1tu);}
\\[-3mm]
\item[] {\tt In [{21}]: ex2 = t\_derivative(E(H,T-t,S(t,u),S1tu),t);}
\\[-3mm]
\item[] {\tt In [{22}]: reduce\_order(expand(ex1-ex2))}
\item[] {\tt Out[{22}]: $0$}
\end{itemize}

\subsection*{Check~\eqref{check3}}
We verify identity~\eqref{check3},
\begin{align*}
&\partial_t\,\nS^{(1)}(t,u) -H'(\nS(t,u))\cdot \nS^{(1)}(t,u) \\
&= \Big[ \partial_2\tilde\nS^{(1)}(t,v)\cdot A(v)
         -A'(\nE_B(t,v))\cdot\tilde\nS^{(1)}(t,v)+[B,A](\nE_B(t,v)) \Big]_{v=\nE_A(t,u)}.
\end{align*}

\begin{itemize}

\item[] {\tt In [23]: {S2tu} {=} {t\_derivative}{(}{S1tu}{,}{t}{)}{-}{A}{(}{Stu}{,}{S1tu}{)}{-}{B}{(}{Stu}{,}{S1tu}{)}}

        \qquad\qquad~\, {\tt {ex1} {=} {S2tu}};
%
\\[-3mm]
\item[] {\tt In [24]: {S2tv} {=} {(differential}{(}{S1tv}{,}{v}{,}{A}{(}{v}{)}{)}{-}{A}{(}{E}{(}{B}{,}{t}{,}{v}{)}{,}{S1tv}{)}

        \qquad\qquad\qquad\qquad~{+}{commutator}{(}{B}{,}{A}{,}{E}{(}{B}{,}{t}{,}{v}{)}{))}

        \qquad\qquad\, {ex2} {=} {substitute}{(}{S2tv}{,}{v}{,}{E}{(}{A}{,}{t}{,}{u}{)}{)};}
%
\\[-3mm]
\item[] {\tt In [{25}]: {expand}{(}{ex1}{-}{ex2}{)}}
\item[] {\tt Out[{25}]: $0$}
\end{itemize}

\subsection*{Check~\eqref{check4}}
We verify identity~\eqref{check4},
\begin{equation*}
\begin{aligned}
\partial_t\tilde\nS^{(2)}(t,v)
&= B'(\nE_B(t,v))\cdot\tilde\nS^{(2)}(t,v)\\
& \quad {} +B''(\nE_B(t,v))(\tilde\nS^{(1)}(t,v),\tilde\nS^{(1)}(t,v))\\
& \quad {} -[B,[B,A]](\nE_B(t,v)) -[A,[B,A]](\nE_B(t,v))\nonumber\\
& \quad {} +2[B,A]'(\nE_B(t,v))\cdot\tilde\nS^{(1)}(t,v).
\end{aligned}
\end{equation*}

\begin{itemize}

\item[] {\tt In [{26}]: {ex1} {=} {(}{B}{(}{E}{(}{B}{,}{t}{,}{v}{)}{,}{S2tv}{)}
              {+} {B}{(}{E}{(}{B}{,}{t}{,}{v}{)}{,}{S1tv}{,}{S1tv}{)}

              \quad\qquad\qquad\qquad~{-}{commutator}{(}{B}{,}{B}{,}{A}{,}{E}{(}{B}{,}{t}{,}{v}{)}{)
              
              \quad\qquad\qquad\qquad~{-}{commutator}{(}{A}{,}{B}{,}{A}{,}{E}{(}{B}{,}{t}{,}{v}{)}{)}}

              \quad\qquad\qquad\qquad~{+}{2}{*}{substitute}{(}{differential}

              \quad\qquad\qquad\qquad~{(}{commutator}{(}{B}{,}{A}{,}{w}{)}{,}{w}{,}{S1tv}{)}{,}{w}{,}{E}{(}{B}{,}{t}{,}{v}{)}{)}{)};}
\\[-3mm]
\item[] {\tt In [{27}]: {ex2} {=} {t\_derivative}{(}{S2tv}{,}{t}{)};}
%
\\[-3mm]
\item[] {\tt In [{28}]: {expand}{(}{ex1}{-}{ex2}{)}}
\item[] {\tt Out[{28}]: $0$}
\end{itemize}

\subsection*{Check~\eqref{check5}}
We verify identity~\eqref{check5},
\begin{equation*}
\begin{aligned}
  &\partial_\tau\Big(
  \partial_2\nE_H(t-\tau,\nS(\tau,u))\cdot
  \partial_2\nE_B(\tau,\nE_A(\tau,u))\cdot [B,A](\nE_A(\tau,u))
  \Big)  =\\
&\qquad \Big[\partial_2\nE_H(t-\tau,\nE_B(\tau,v))\cdot\partial_2\nE_B(\tau,v)\cdot[[B,A],A](v)\\
& \qquad\qquad {} + \partial_2\nE_H(t-\tau,\nE_B(\tau,v))\cdot\partial_2\tilde\nS^{(1)}(\tau,v)\cdot[B,A](v)\\
& \qquad\qquad {} + \partial_2^2\nE_H(t-\tau,\nE_B(\tau,v))\big(\tilde\nS^{(1)}(\tau,v), \partial_2\nE_B(\tau,v)\cdot[B,A](v)\big)
\Big]_{v=\nE_A(t,u)}.
\end{aligned}
\end{equation*}

\begin{itemize}

\item[] {\tt In [{29}]: {ex1} {=} {t\_derivative}{(}{E}{(}{H}{,}{T}{-}{t}{,}{Stu}{,}{E}{(}{B}{,}{t}{,}{E}{(}{A}{,}{t}{,}{u}{)}{,}

\quad\qquad\qquad\qquad~{commutator}{(}{B}{,}{A}{,}{E}{(}{A}{,}{t}{,}{u}{)}{)}{)}{)}{,}{t}{)};}
%
\item[] {\tt In [{30}]: {ex2} {=} {(substitute}{(}{-}{E}{(}{H}{,}{T}{-}{t}{,}{E}{(}{B}{,}{t}{,}{v}{)}{,}{E}{(}{B}{,}{t}{,}{v}{,}

\quad\qquad\qquad\qquad~{commutator}{(}{A}{,}{B}{,}{A}{,}{v}{)}{)}{)}{+}{E}{(}{H}{,}{T}{-}{t}{,}{E}{(}{B}{,}{t}{,}{v}{)}{,}

\quad\qquad\qquad\qquad~{differential}{(}{S1tv}{,}{v}{,}{commutator}{(}{B}{,}{A}{,}{v}{)}{)}{)}

\quad\qquad\qquad\qquad~{+}{E}{(}{H}{,}{T}{-}{t}{,}{E}{(}{B}{,}{t}{,}{v}{)}{,}{S1tv}{,}{E}{(}{B}{,}{t}{,}{v}{,}
                         
\quad\qquad\qquad\qquad~{commutator}{(}{B}{,}{A}{,}{v}{)}{)}{)}{,}{v}{,}{E}{(}{A}{,}{t}{,}{u}{)}{))};}
\\[-3mm]
\item[] {\tt In [{31}]: {diff} {=} {ex1}{-}{ex2}

        \qquad\qquad~{diff} {=} {FE2DEF}{(}diff{)}

        \qquad\qquad~{diff} {=} {substitute}{(}{diff}{,}{H}{,}{A}{(}{v}{)}{+}{B}{(}{v}{)}{,}{v}{)}

        \qquad\qquad~{diff} {=} {expand}{(}{reduce\_order}{(}{diff}{)}{)}}
\item[] {\tt Out[{31}]: $0$}
\end{itemize}

\section{Elementary differentials}

To further demonstrate the functionality and correctness of our package,
we consider the
elementary differentials obtained by repeated  differentiation of a differential equation,
\begin{align*}
y'(t)&=F(y(t)), \\
y''(t)&=F'(y(t))\cdot F(y(t)), \\
y'''(t)&=F''(y(t))(F(y(t),F(y(t)) + F'(y(t))\cdot F'(y(t))\cdot F(y(t)), \\
&~~\vdots&
\end{align*}
The number of terms in these expressions are available from the literature, see~\cite[Table~2.1]{haireretal93_2nd_ed}
or~\cite{IntSeq}.
The following notebook generates the elementary differentials  and counts their number:
\begin{itemize}
\item[] {\tt In [{1}]: {using} {Flows}}
\\[-3mm]
\item[] {\tt In [{2}]: {@}{t\_vars} {t};

        \qquad\qquad {@}{x\_vars} {u};
        
        \qquad\qquad {@}{funs} {F}};
\\[-3mm]
\item[] {\tt In [{3}]: {ex} {=} {E}{(}{F}{,}{t}{,}{u}{)}}
\item[] {\tt Out[{3}]: $\mathcal{E}_{F}(t,u)$}
\\[-3mm]
\item[] {\tt In [{4}]: {ex} {=} {t\_derivative}{(}{ex}{,}{t}{)}}
\item[] {\tt Out[{4}]: $F(\mathcal{E}_{F}(t,u))$}
\\[-3mm]
\item[] {\tt In [{5}]: {ex} {=} {t\_derivative}{(}{ex}{,}{t}{)}}
\item[] {\tt Out[{5}]:}
    \(\;F'(\mathcal{E}_{F}(t,u))\cdot F(\mathcal{E}_{F}(t,u))\)
\\[-3mm]
\item[] {\tt In [{6}]: {ex} {=} {t\_derivative}{(}{ex}{,}{t}{)}}
\item[] {\tt Out[{6}]:}
\vspace{-7.25mm}
\begin{align*}
&F'(\mathcal{E}_{F}(t,u))\cdot F'(\mathcal{E}_{F}(t,u))\cdot F(\mathcal{E}_{F}(t,u)) \\[-0.5\jot]
&+F''(\mathcal{E}_{F}(t,u))(F(\mathcal{E}_{F}(t,u)),F(\mathcal{E}_{F}(t,u)))
\qquad\qquad\qquad~~
\end{align*}
\pagebreak
\item[] {\tt In [{7}]: {ex} {=} {t\_derivative}{(}{ex}{,}{t}{)}}
\item[] {\tt Out[{7}]:}
\vspace{-7.25mm}
\begin{align*}
    &3F''(\mathcal{E}_{F}(t,u))(F'(\mathcal{E}_{F}(t,u))\cdot F(\mathcal{E}_{F}(t,u)),F(\mathcal{E}_{F}(t,u)))\quad\,\,\,\,\\[-0.5\jot]
    &+F'''(\mathcal{E}_{F}(t,u))(F(\mathcal{E}_{F}(t,u)),F(\mathcal{E}_{F}(t,u)),F(\mathcal{E}_{F}(t,u)))\\[-0.5\jot]
    &+F'(\mathcal{E}_{F}(t,u))\cdot (F'(\mathcal{E}_{F}(t,u))\cdot F'(\mathcal{E}_{F}(t,u))\cdot F(\mathcal{E}_{F}(t,u))\\[-0.5\jot]
    &+F''(\mathcal{E}_{F}(t,u))(F(\mathcal{E}_{F}(t,u)),F(\mathcal{E}_{F}(t,u))))
\end{align*}
\\[-11mm]
\end{itemize}
This is not yet fully expanded.
It is a linear combination consisting of 3 terms:
\begin{itemize}
\item[] {\tt In [{8}]: {length}{(}{ex}{.}{terms}{)}}
\item[] {\tt Out[{8}]: 3}
\end{itemize}
If we expand it, we obtain a linear combination
of~4 terms, corresponding to the 4~elementary differentials
(Butcher trees) of order~4:
\begin{itemize}
\item[] {\tt In [{9}]:   \hspace{-0.3mm} {ex} {=} {expand}{(}{ex}{)}}
\item[] {\tt Out[{9}]:}
\vspace{-7.25mm}
    \begin{align*}
    &3F''(\mathcal{E}_{F}(t,u))(F'(\mathcal{E}_{F}(t,u))\cdot F(\mathcal{E}_{F}(t,u)),F(\mathcal{E}_{F}(t,u))) \\[-0.5\jot]
    &+F'(\mathcal{E}_{F}(t,u))\cdot F'(\mathcal{E}_{F}(t,u))\cdot F'(\mathcal{E}_{F}(t,u))\cdot F(\mathcal{E}_{F}(t,u)) \\[-0.5\jot]
    &+F'''(\mathcal{E}_{F}(t,u))(F(\mathcal{E}_{F}(t,u)),F(\mathcal{E}_{F}(t,u)),F(\mathcal{E}_{F}(t,u))) \\[-0.5\jot]
    &+F'(\mathcal{E}_{F}(t,u))\cdot F''(\mathcal{E}_{F}(t,u))(F(\mathcal{E}_{F}(t,u)),F(\mathcal{E}_{F}(t,u)))
    \end{align*}
\\[-11mm]
\item[] {\tt In [{10}]: {length}{(}{ex}{.}{terms}{)}}
\item[] {\tt Out[{10}]: 4}
\\[-3mm]
\item[] {\tt In [{11}]: {ex} {=} {expand}{(}{t\_derivative}{(}{ex}{,}{t}{)}{)};}
\\[-3mm]
\item[] {\tt In [{12}]: {length}{(}{ex}{.}{terms}{)}}
\item[] {\tt Out[{12}]: 9}
\\[-3mm]
\item[] {\tt In [{13}]: {ex} {=} {expand}{(}{t\_derivative}{(}{ex}{,}{t}{)}{)};}
\\[-3mm]
\item[] {\tt In [{14}]: {length}{(}{ex}{.}{terms}{)}}
\item[] {\tt Out[{14}]: 20}
\\[-3mm]
\item[] {\tt In [{15}]: {ex} {=} {expand}{(}{t\_derivative}{(}{ex}{,}{t}{)}{);}}
\\[-3mm]
\item[] {\tt In [{16}]: {length}{(}{ex}{.}{terms}{)}}
\item[] {\tt Out[{16}]: 48}
\vspace{-2.0mm}
\end{itemize}
\begin{verbatim}
   In [17]: ex = E(F,t,u)
            ex = expand(t_derivative(ex,t))
            println("order\t#terms")
            println("---------------")
            println(1,"\t",1)
            ex = expand(t_derivative(ex,t))
            println(2,"\t",1)
            for k=3:16
                ex = expand(t_derivative(ex,t))
                println(k,"\t",length(ex.terms))
            end
\end{verbatim}
\newpage
\begin{verbatim}
     order   #terms
     ---------------
     1       1
     2       1
     3       2
     4       4
     5       9
     6       20
     7       48
     8       115
     9       286
     10      719
     11      1842
     12      4766
     13      12486
     14      32973
     15      87811
     16      235381
\end{verbatim}

\section{Acknowledgments}

This work was supported in part by the projects P24157-N13 of
the Austrian Science Fund (FWF) and MA14-002 of the
Vienna Science and Technology Fund (WWTF).

\bibliographystyle{splncs03}
\bibliography{NEW,schroedinger,num,gen,books}

\end{document}